\documentclass[10pt]{amsart}
\usepackage[all]{xy}
\usepackage{amsmath}
\usepackage{amsfonts}
\usepackage{amssymb}
\usepackage{amscd}
\usepackage{verbatim}
\newtheorem{theorem}{Theorem} [section]
\newtheorem{lemma} [theorem] {Lemma}
\newtheorem{proposition} [theorem] {Proposition}

\newtheorem{corollary} [theorem] {Corollary}

\begin{document}
\title{Hartogs extension theorems on Stein spaces}
\author{Nils \O vrelid and Sophia Vassiliadou}
\thanks{{\em 2000 Mathematics Subject Classification:} 32B10, 32J25, 32W05, 14C30}
\thanks{The research of the second author is partially supported by NSF grant DMS-0712795}
\keywords{Cauchy-Riemann equation, Singularity, Cohomology groups}
\address{Dept. of Mathematics\\University of Oslo\\
P.B 1053 Blindern, Oslo, N-0316 NORWAY}
\address{Dept. of Mathematics\\Georgetown University\\
Washington, DC 20057 USA} \email{nilsov@math.uio.no,\;
sv46@georgetown.edu}
\begin{abstract} \noindent We discuss various known generalizations of the
classical Hartogs  extension theorem on Stein spaces with arbitrary
singularities and present an analytic proof based on
$\overline\partial$-methods.
\end{abstract}

\maketitle
\medskip
\noindent \section{Introduction}

\medskip
\noindent\noindent
\medskip
\noindent The classical Hartogs extension theorem  is usually stated
as follows:

\begin{theorem}(Hartogs)  Let $\Omega$ be a domain in $\mathbb{C}^n$ ($n>1$) and let $K$ be a compact subset
of $\Omega$ such that $\Omega\setminus K$ is connected. Then  any
holomorphic function $f$ on $\Omega\setminus K$ has a unique
holomorphic extension on $\Omega$.
\end{theorem}

\medskip
\noindent The first analytic proofs of this theorem for
$\Omega\subset\subset \mathbb{C}^n$ used the theory of integral
kernels. The shortest proof so far of Theorem 1.1 is due to
Ehrenpreis and relies on solvability with compact support of the
Cauchy-Riemann equations for  $\overline\partial$-closed smooth
compactly supported $(0,1)$-forms defined in $\mathbb{C}^n$.
Vanishing results of this type are instrumental in obtaining Bochner
type extension theorems as well. Complex manifolds $X$ that possess
this property i.e. $H^{1}_{c}(X, \mathcal{O})=0$ are often said to
have the Bochner-Hartogs property. Examples of such manifolds are
$(n-1)$-complete complex manifolds of dimension $n>1$ (this was
shown by Andreotti-Grauert in \cite{AG} and Andreotti and Vesentini
in \cite{AV}), and $(n-1)$-strictly-hyperconvex K\"ahler manifolds
(proven  by Grauert and Riemenschneider in \cite{GR}). Various
possible Hartogs extension theorems have appeared in the literature
in which one can replace $\mathbb{C}^n$ by an $n$-dimensional Stein
manifold $X$ or, replace the holomorphic sections by sections of a
coherent sheaf on $X$. Similarly, various Bochner extension theorems
have appeared in which one can replace $\mathbb{C}^n$ by a Stein or
$(n-1)$-complete  $n$-dimensional manifold (see section 3.4 in
\cite{AH}). A weak form of Hartogs  extension theorems for singular
spaces first appeared in the seminal 1962 paper of Andreotti and
Grauert (Th$\acute e$or$\grave e$me 15 in \cite{AG}). In 1963, Rossi
proved the following version of a Hartogs  extension theorem for
normal Stein spaces using Bishop's special analytic polyhedra:

\begin{theorem} (Theorem 6.6 in \cite{R}) Let $X$ be a connected, normal Stein space of dimension $n\ge 2$.
Let $K$ be a compact subset of $X$ and let $h$ be a holomorphic
function on $X\setminus K$. There is a unique $H$ holomorphic on $X$
such that $H=h$ on the unbounded component of $X\setminus K$.
\end{theorem}

\medskip
\noindent Later on Laufer,\;\,Harvey, B\v{a}nic\v{a} and
St\v{a}n\v{a}\c{s}il\v{a} generalized Rossi's result. Harvey showed
in \cite{Harv} the following

\begin{theorem} (Theorem 2.13 in \cite{Harv}) Let $X$ be a subvariety in $\mathbb{C}^n$. Suppose $K$ is a compact subset of
$X$ such that $X\setminus K$ has no branches which are contained in
a compact set. If the homological codimension $\text{codh}_x
(\widetilde{\mathcal{O}}_{X})\ge 2$ for all $x\in X$, then the
restriction map

$$
\Gamma(U, \mathcal{O}_X)\to \Gamma (U\setminus K, \mathcal{O}_X)$$

\smallskip
\noindent is bijective for all neighborhoods $U$ of $K$ in $X$.
\newline
Conversely, if \;$\Gamma(X\setminus K, \mathcal{O}_X)\cong \Gamma(X,
\mathcal{O}_X)$ for a compact holomorphically convex subset of
$\mathbb{C}^n$ with $K\subset X$, then $\text{codh}_x
(\widetilde{\mathcal{O}}_{X})\ge 2$ for all $x\in K$. Here
$\widetilde{\mathcal{O}}_{X}:=i_{*} (\mathcal{O}_X)$ is,  the
trivial extension of $\mathcal{O}_X$ by zero to $\mathbb{C}^n$, and
$i: X\hookrightarrow \mathbb{C}^n$ is the closed immersion.
\end{theorem}

\medskip
\noindent On another note in B\v{a}nic\v{a} and
St\v{a}n\v{a}\c{s}il\v{a}'s book \cite{BS} the following
generalization of Theorem 1.3 is given.

\begin{theorem} (Chapter I, Corollary 4.4 in \cite{BS}) Let $(X, \mathcal{O})$ be a reduced Stein space and $K$ be a
compact set such that $X\setminus K$ has no relatively compact
irreducible components (branches) in $X$. If the homological
codimension $\text{codh}\, (\mathcal{O}_X)\ge
2\footnotemark\footnotetext{Corollary 4.4 on page 42 in \cite{BS} is
stated under a weaker assumption on the homological codimension of
the structure sheaf; namely that $\text{codh}(\mathcal{O}_x)\ge 2$
for all $x\in K$. However a stronger assumption is needed for the
proof of Corollary 4.4 in \cite{BS} to go through. One needs at
least that $\text{codh}\,(\mathcal{O}_x)\ge 2$ for all $x\in U$
where $U$ is  a Stein neighborhood of $K$.}$, then the restriction
map

$$\Gamma(X, \mathcal{O})\to \Gamma(X\setminus K, \mathcal{O})$$ is bijective.

\noindent Conversely, if $X$ is Stein and $K$ a holomorphically
convex compact such that the map $\Gamma(X, \mathcal{O})\to
\Gamma(X\setminus K, \mathcal{O})$ is bijective, then
$\text{codh}\,(\mathcal{O}_x)\ge 2$ for all $x\in K$.
\end{theorem}

\medskip
\noindent In 2007, Merker and Porten established a version of a
Hartogs  extension theorem for $(n-1)$-complete (possibly singular)
complex spaces using the method of pushing disks and Morse theory.
More precisely they proved

\begin{theorem} (Theorem 2.2 in \cite{MP}) Let $X$ be a connected $(n-1)$-complete normal space
of dimension $n\ge 2$. Then, for every domain $\Omega\subset X$ and
every compact set $K\subset \Omega$ with $\Omega\setminus K$
connected, holomorphic functions on $\Omega\setminus K$ extend
holomorphically and uniquely to $\Omega$.
\end{theorem}

\medskip
\noindent Ruppenthal proved in \cite{Rup1} Theorem 1.5 for
1-complete (Stein) spaces with isolated singularities, using the
classical $\overline\partial$-method of Ehrenpreis, ideas from an
earlier paper of ours  \cite{FOV1} combined with Serre duality.

\medskip
\noindent In the case of Stein spaces with singularities it is
natural to examine the relationship of the above versions of the
Hartogs  extension theorem and also try to generalize Ehrenpreis
approach for Stein spaces with arbitrary singularities. We shall
show that Theorem 1.4 and Theorem 1.5 for Stein spaces are in fact
equivalent; we will further show that Theorem 1.5 for Stein spaces
implies the version of Theorem 1.4 under the {\it{weaker}}
assumption that $\text{codh}\,(\mathcal{O}_x)\ge 2$ for all $x\in
K$. We shall also present in this paper an analytic approach to
obtaining Hartogs extension theorems based on our earlier work. It
may be worthwhile to state the following corollary which will tell
us that only appropriate neighborhoods of $K$ need to be considered
for extension theorems to hold:

\begin{corollary} Let $(X, \mathcal{O})$ be a reduced complex
analytic space $K$ be a compact subset of $X$ such that
$\text{codh}(\mathcal{O}_{X,x})\ge 2$ for all $x\in K$. Let $D$ be
an open neighborhood of $K$. If $K$ has a Stein neighborhood $U$ in
$X$ such that $U\setminus K$ has no irreducible components that are
relatively compact in $U$ then the restriction map

$$
\Gamma (D, \mathcal{O})\to \Gamma(D\setminus K, \mathcal{O})
$$

\noindent is bijective.
\end{corollary}

\medskip
\noindent The organization of the paper is as follows. In section 2,
we will give an overview of the ideas used to prove Theorems 1.3,
1.4. In section 3 we will prove the equivalence of the various
Hartogs  extension theorems, the implication that Theorem 1.5
implies Corollary 4.4  as stated in \cite{BS} and show that the
conditions imposed in Theorem 1.4 are necessary. In section 4 we
will give a different proof of the Hartogs  extension for Stein
spaces $X$ that satisfy $\text{codh}(\mathcal{O}_X)\ge 2$ based on
Andreotti-Grauert's work. Section 5 will streamline the
$\overline\partial$-approach.

\medskip
\noindent{\it{Remark:}} This paper was completed in September of
2008. Around the same time Andersson and Samuelsson obtained in
\cite{AS} (Theorem 1.4) a $\overline\partial$-proof of a Hartogs
extension theorem on Stein spaces with arbitrary singularities,
using residue theory calculus. In November of 2008, Coltoiu \cite{C}
and Ruppenthal \cite{Rup2} independently obtained a
$\overline\partial$-theoretic proof of a Hartogs extension theorem
on cohomological $(n-1)$-complete (resp. $(n-1)$-complete) spaces.
The key ingredient in the proof is a vanishing of the higher direct
images of the sheaf of canonical forms of an appropriate
desingularization $\tilde{X}$ of the $(n-1)$-complete complex space
$X$. This subtle vanishing result which was obtained by Takegoshi
\cite{T}, easily yields the vanishing of $H^1_c( \tilde{X},
\mathcal{O})$ which is needed for the Ehrenpreis method to carry
over. Our analytic approach is based on more general weighted
$L^2$-solvability results for $\overline\partial$-closed, compactly
supported forms (Theorem 5.3) which are of independent interest. In
a forthcoming work \cite{OV} we obtain semiglobal results for
$\overline\partial$ on $q$-complete spaces and as a corollary we
obtain another proof of a Hartogs extension theorem on
$(n-1)$-complete spaces.

\section{Preliminaries}

\subsection{Preliminaries on cohomology with support}

\medskip
\noindent Let us recall the obstructions to extending holomorphic
sections from $X\setminus K$ to $X$. For a sheaf of abelian groups
$\mathcal{S}$ on $X$, let $\Gamma_K (X,\mathcal{S})$ denote the
sections on $X$ with support in $K$. Consider a flabby resolution of
$\mathcal{O};\;\;\; 0\to \mathcal{O}\to
\mathcal{C}^{0}\overset{d_0}\to \mathcal{C}^{1}\overset{d_1}\to
\cdots$. The cohomology groups with support in $K$ are defined by
$H^{*}_K (X, \mathcal{O}):=H^{*} (\Gamma_K \left(X,
\mathcal{C}^{\bullet})\right)$, i.e. they are the cohomology groups
of the complex $\left(\Gamma_K(X, \mathcal{C}^{k}), d_k \right)$.
Since each $\mathcal{C}^{\bullet}$ is flabby, we have a short exact
sequence $0\to \Gamma_K (X, \mathcal{C}^{\bullet})\to \Gamma (X,
\mathcal{C}^{\bullet})\to \Gamma (X\setminus K,
\mathcal{C}^{\bullet})\to 0$. This induces a long exact sequence on
cohomology

\begin{equation}\label{eq:oles}
0\to H^0_K(X,\mathcal{O})\to H^0(X, \mathcal{O})\to H^0(X\setminus
K, \mathcal{O})\to H^1_K(X,\mathcal{O})\to H^1(X,\mathcal{O})\to...
\end{equation}

\noindent Since $X$ is a Stein subvariety $H^i(X,\mathcal{O})=\{0\}$
for all $i\ge 1$. Hence the above sequence becomes

\begin{equation}\label{eq:les}
0\to H^0_K(X,\mathcal{O})\to H^0(X, \mathcal{O})\to H^0(X\setminus
K, \mathcal{O})\to H^1_K(X,\mathcal{O})\to 0
\end{equation}

\noindent It is clear from the above that the obstructions to
extending holomorphic sections from $X\setminus K$ to $X$ come from
nontrivial $H^i_K(X,\mathcal{O})$ for $i=0,\,1$.

\medskip
\noindent It is a standard fact from sheaf cohomology theory that
$H^i_K(X,\mathcal{O})\cong H^i_K(U, \mathcal{O})$ where $U$ is an
open neighborhood of $K$ in $X$. The fact that
$H^i_K(U,\mathcal{O})$ is independent of the neighborhood $U$ of $K$
is referred to as excision.

\medskip
\noindent Later on we shall need a more general version of
$(\ref{eq:oles})$. Let $L$ be a closed neighborhood of $K$ in $X$.
Then we have the following exact sequence:

\begin{equation}\label{eq:red2}
0\to H^0_K(X,\mathcal{O})\to H^0_L(X, \mathcal{O})\to
H^0_{L\setminus K}(X\setminus K, \mathcal{O})\to H^1_K(X,
\mathcal{O})\to H^1_{L}(X,\mathcal{O})\to.. \end{equation}

\subsection{Preliminaries on homological codimension}

\noindent Let us recall the  definition  and some properties of the
homological codimension or profundity of a coherent sheaf. First of
all,

\medskip
\noindent $\bullet_1$  On a noetherian ring $R$ with an $R$-module
$M$ of finite type and $\mathcal{\alpha}$ an ideal, a sequence of
elements $f_1, \cdots, f_q\in \mathcal{\alpha}$ is called a regular
$M$-sequence if $f_i$ is not a zero divisor of $M/{\sum_{j=1}^{i-1}
f_j\, M}$ (for $i=1$, it means that $f_1$ is not a zero divisor of
$M$). The maximum length of a regular $M$-sequences is denoted by
$\text{prof}_{\mathcal{\alpha}} M$. If $R$ is a local ring with
maximal ideal ${m}$ then $\text{prof}_{m} M$ is denoted by
$\text{codh}_R M$ and is called the homological codimension of $M$.
An equivalent way to describe the homological codimension of a
module $M$ of a local noetherian ring  is the following:
$\text{codh}_R M=\text{inf} \{i;\;\;\text{Ext}^i (k, \, M)\neq 0\}$
where $k=R/m$.

\smallskip
\noindent For a general complex space $(X,\mathcal{O})$ and a
coherent analytic sheaf $\mathcal{F}$ on $X$, we define
$\text{codh}\,\mathcal{F}=\text{inf}_{x\in X} \text{codh}_x
\mathcal{F}$ where  $\text{codh}_x
\mathcal{F}=\text{codh}_{\mathcal{O}_x} \mathcal{F}_x$, if
$\mathcal{F}_x\neq 0$ and $\infty$ otherwise.

\medskip
\noindent $\bullet_2$ If $i:X\rightarrow Y$ is a closed immersion of
complex spaces and $\mathcal{F}$ is a coherent sheaf on $X$ then
$i_{*}(\mathcal{F})$, the trivial extension of $\mathcal{F}$ by zero
outside $X$, is a coherent sheaf in $Y$ and we have
$\text{codh}(\mathcal{F}_x)=\text{codh}\,(i_{*}(\mathcal{F})_{i(x)})$
for all $x\in X$.


\medskip
\noindent Hartogs extension theorems for some singular Stein spaces
follow as a special case of more general vanishing theorems of
cohomology with supports. In what follows we outline the basic ideas
behind these vanishing theorems.

\subsection{Vanishing on cohomology with supports and homological codimension bounds}

\medskip
\noindent The aim in this section is to recall the relationship
between homological codimension bounds and vanishing of cohomology
with supports.

\begin{theorem} (Corollary 3.1, page 37 in \cite{BS}) If $X$ is a Stein space, $\mathcal{F}$ is a
coherent sheaf on $X$, $K$ is a holomorphically convex compact and
$r\ge 0$ an integer. Then $\text{codh}\,(\mathcal{F}_x)\ge r+1$ for
all $x\in K$ if and only if $H^{i}_K (X, \,\mathcal{F})=\{0\}$ for
all $i\le r$, where $H^i_K$ are the $i$-th cohomology groups with
support on $K$.
\end{theorem}

\medskip
\noindent Since $K$ is a holomorphically convex compact in a Stein
space, it has a neighborhood basis of Oka-Weil domains i.e. $K$ has
a Stein neighborhood that can be thought of as a holomorphic
subvariety of a polydisk in some numerical affine space. So without
loss of generality we can assume that $i: X\hookrightarrow
\mathbb{C}^n$ is a closed immersion. Let
$\mathcal{G}:=i_*(\mathcal{F})$, the trivial extension by zero
outside $X$. The sheaf $\mathcal{G}$ is  coherent. If $K$ is
holomorphically convex subset of $X$ then $L:=i(K)$ is
holomorphically convex compact in $\mathbb{C}^n$ and
$H^i_K(X,\mathcal{F})\cong H^i_L(\mathbb{C}^n,\, \mathcal{G})$.
Hence the proof of Theorem 2.1 reduces to the case of Stein
manifolds. So what we really want to show is that on the
$n$-dimensional Stein manifold $X:=\mathbb{C}^n$ with $\mathcal{F}$
a coherent sheaf and $K$ holomorphically convex compact,

$$\text{codh}(\mathcal{F}_x)\ge r+1\;\;\text{for\;\; all}\;\; x\in K\;\;\Leftrightarrow\;\;H^i_K(X,\,\mathcal{F})=0\;\;\text{for\;\; all}\;\; i\le r.$$

\medskip
\noindent The standard proof relies into two main ingredients. First
one shows that the homological codimension bounds imply the
vanishing of certain extension sheaves. Then one proves that the
cohomology groups $H^i_K(X,\,\mathcal{F})$ with their natural
topology are the strong duals of the spaces $H^0 (K,
{\mathcal{E}xt}^{n-i}_{\mathcal{O}}(\mathcal{F},\,\mathcal{O})\otimes_{\mathcal{O}}
\Omega^n)$, where $\Omega^n$ is the sheaf of holomorphic $(n,0)$
forms on $X$.

\medskip
\noindent {\bf{Claim I:}} $\text{codh}(\mathcal{F}_x)\ge r+1$ for
all $x\in K \Leftrightarrow
{\mathcal{E}xt}^{n-i}_{\mathcal{O}}(\mathcal{F}, \mathcal{O})=0$
\;over $K$ and for all $i\le r$.

\smallskip
\noindent Since $K$ is a holomorphically convex compact in
$\mathbb{C}^n$ we can find an open Stein neighborhood $U$ of $K$
over which there exists an exact sequence of holomorphic free
sheaves of the form

\begin{equation}\label{eq:proj}
0\longrightarrow
\mathcal{O}^{\nu_r}\overset{\delta_r}\longrightarrow
\mathcal{O}^{\nu_{r-1}}\longrightarrow \cdots\longrightarrow
\mathcal{O}^{\nu_{0}} \overset{\delta_0}\longrightarrow
\mathcal{F}_{|_{U}}\longrightarrow 0.
\end{equation}

\noindent We can always choose this sequence to be of minimal
length. Let
$\check{\mathcal{F}}:=\text{Hom}_{\mathcal{O}}(\mathcal{F},
\mathcal{O})$ be the dual of $\mathcal{F}$. Since
$\text{Hom}_{\mathcal{O}}(\cdot,\; \mathcal{O})$ is a contravariant
functor by applying it to $(\ref{eq:proj})$ we obtain the following
sequence over $U$

\begin{equation}\label{eq:defext}
0\longrightarrow
\check{\mathcal{F}}\overset{\delta^t_0}\longrightarrow
\mathcal{O}^{\nu_{0}}\overset{\delta^t_{1}}\longrightarrow
\mathcal{O}^{\nu_1}\longrightarrow \cdots
\overset{\delta^t_{r}}\longrightarrow
\mathcal{O}^{\nu_r}\longrightarrow 0.
\end{equation}

\noindent The above complex is not exact-except at
$\check{\mathcal{F}}$ and $\mathcal{O}_{\nu_0}$. Let us define
$\mathcal{E}xt^i_{\mathcal{O}}(\mathcal{F},\,\mathcal{O}):=\text{kern}
(\delta^t_{i+1})/\text{Im}(\delta^t_i)$ for $i>0$. These sheaves
measure the degree to which the above complex fails to be exact.
They are coherent sheaves on $K$. Let
$\text{codh}_K(\mathcal{F}):=\text{inf}_{x\in
K}\,\text{codh}(\mathcal{F}_x)$. By the assumption of the claim
$\text{codh}_K (\mathcal{F})\ge r+1$. Then $n-\text{codh}_K
(\mathcal{F})\le n-r-1<n-r$. Hence for all $i\le
r,\;\;{\mathcal{E}xt}^{n-i}_{\mathcal{O}}(\mathcal{F},\,\mathcal{O})=0$
over $K$.

\medskip
\noindent {\bf{Claim II:}} $H^i_K(X,\mathcal{F})$ is the strong dual
of $H^0\left(K,\, \mathcal{E}xt^{n-i}_{\mathcal{O}}(\mathcal{F},
\,\mathcal{O})\otimes \Omega^n\right)$.

\medskip
\noindent The proof of the above claim proceeds in several steps.
First you consider the case of locally free sheaves $\mathcal{F}$
and you show that i) $H^i_K(X, \mathcal{F})=0$ if $i\neq n$,\; ii)
$H^0(K,\,\check{\mathcal{F}})$ with the natural inductive limit
topology is a (DFS) space (strong dual of a Fr\'echet space ) and
iii)  $H^n_K(X, \mathcal{F})$ with its natural topology is the
strong dual of $H^0 (K, \check{\mathcal{F}}\otimes \Omega^n)$ (see
theorem 5.9 in \cite{ST}). Serre's duality theorem plays a critical
role in the proofs of the above statements. For a general coherent
sheaf $\mathcal{F}$ one uses a short exact sequence of the form
$0\to \mathcal{G}\to \mathcal{O}^{s_0}\to \mathcal{F}\to 0$ over a
Stein neighborhood of $K$ (such a sequence always exist since
$\mathcal{F}$ is coherent) and apply induction on the homological
codimension of $\mathcal{F}$ over $K$ to prove that $H^0\left(K,
\mathcal{E}xt^{n-i}_{\mathcal{O}}(\mathcal{F},\,\mathcal{O})\otimes
\Omega^n\right)$ with its natural inductive limit topology is a
(DFS) space and then show that $H^i_K (X,\mathcal{F})$ with its
natural topology is the strong dual of $H^0\left(K,
\mathcal{E}xt^{n-i}_ {\mathcal{O}}(\mathcal{F},\,\mathcal{O})\otimes
\Omega^n\right)$ (see theorem 5.12 in \cite{ST}).

\medskip
\noindent{\bf{Historical remarks:}} (i)-(iii) in the previous
paragraph are special cases of more general theorems proven by
Martineau in  \cite{M}. The passage from locally free sheaves to
coherent sheaves is due to Harvey \cite{Harv}.

\medskip
\noindent Theorem 2.1 follows immediately from Claims I and II.

\subsection{Proof of Theorems 1.3, 1.4}
\medskip
\noindent  {\it{Reduction steps.}} $\alpha$)  Let us suppose that
the first half of Theorem 1.3 is true for $U=X$. This implies that
$H^0_K(X,\mathcal{O})=H^1_K(X, \mathcal{O})=\{0\}$. Now, let $U$ be
any open neighborhood of $K$. Consider the long exact sequence

$$
0\to H^0_K(U, \mathcal{O})\to H^0(U, \mathcal{O})\to H^0 (U\setminus
K, \mathcal{O})\to H^1_K(U, \mathcal{O})\to...
$$

\medskip
\noindent From the excision property we know that $H^i_K(X,
\mathcal{O})\cong H^i_K(U, \mathcal{O})$ for any open neighborhood
$U$ of $K$. Since the former cohomology groups vanish these will
imply that $H^i_K(U, \mathcal{O})=\{0\}$ for $i=0,\,1$. Hence,
$H^0(U, \mathcal{O})\cong H^0(U\setminus K, \mathcal{O})$.

\medskip
\noindent$\beta)$ We need to show now that the theorem is true for
$U=X$. If $K$ is a holomorphically convex compact in $X$ such that
$\text{codh}(\mathcal{O}_x)\ge 2$ for all $x\in K$, then due to
Theorem 2.1 we have that $\Gamma(X\setminus K, \mathcal{O})\cong
\Gamma(X,\mathcal{O})$ or equivalently that $H^i_K(X,
\mathcal{O})=\{0\}$ for $i=0,\,1$.

\medskip
\noindent Now, let $K$ be as in Theorem 1.2. Apply Theorem 2.1 to
its holomorphically convex hull $\hat{K}_X$ in $X$. This will imply
that $H^i_{\hat{K}_X}(X,\,\mathcal{O})=0$ for $i=0,\,1$. Using these
vanishing results and $(\ref{eq:red2})$ we deduce that
$H^0_K(X,\mathcal{O})=\{0\}$ and $H^1_K(X,\mathcal{O})\cong
H^0_{{\hat{K}_X}\setminus K}(X\setminus K, \mathcal{O})$.  The
assumption that $X\setminus K$ has no relatively compact branches in
$X$ implies that $H^0_{{\hat{K}_X}\setminus K}(X\setminus K,
\mathcal{O})=\{0\}$ and hence $H^1_K(X, \mathcal{O})=\{0\}$.

\medskip
\noindent The prove the converse conclusion we refer to Theorem 2.1
and $(\ref{eq:les})$.

\medskip
\noindent The proof of Theorem 1.4 will proceed along the same
lines.


\section{Equivalence of Theorems 1.4 and 1.5 for Stein spaces}

\subsection{Preliminaries on topological properties of complements of compact subsets
on connected normal Stein spaces}

\medskip
\noindent Let $(X, \mathcal{O})$ be a connected normal (thus
irreducible) Stein space of dimension $n\ge 2$.  Let $K$ be a
compact subset of $X$ and let $\hat{K}$ denote the holomorphically
convex hull of $K$ in $X$. Let $\{V_i\}_{i\in I}$ (resp.
$\{U_j\}_{j\in J}$) be the connected components of $X\setminus K$
that are relatively compact (resp. unbounded) in $X$. Then
$X\setminus K=(\cup_{i\in I} V_i) \bigcup (\cup_{j\in J} U_j)$.

\medskip
\noindent {\bf{Claim:}} Let $K$ be a holomorphically convex compact
subset of a connected normal Stein space $X$ of dimension $n\ge 2$.
Then, $X\setminus K$ has no relatively compact connected components.

\medskip
\noindent{\it{Proof of Claim:}} Suppose on the contrary that
$X\setminus K$ has a relatively compact connected component $V$.
Now, $X\setminus K$ is open thus normal and  $V$ is open, connected
with $\overline{V}$  compact. Any $f\in \Gamma(X, \mathcal{O})$ is
continuous on $\overline{V}$, hence it attains its maximum on $p\in
\overline{V}$. If $p\in V$ then we would be able to find a small,
Stein neighborhood $W$ of $p$ such that $f$ attains a local maximum
there. Hence $f$ would be constant near $p$ (by the local maximum
principle, page 109 in \cite{GRem}). Since $X$ is irreducible  by
the Identit\"atsatz (page 170 in \cite{GRem}) $f$ should be constant
on $X$. Hence $f$ attains its maximum on $\partial V$. But $\partial
V \subset K$ and for all $z\in V$ and all $f\in \Gamma(X,
\mathcal{O})$ we have: $|f(z)|\le |f(p)|\le \text{sup}_K |f|$. This
last inequality implies that $V\subset \hat{K}=K$ since $K$ is
holomorphically convex compact and hence $V=\emptyset$. So
$X\setminus K$ has no relatively compact connected components.

\medskip
\noindent{\bf{Remark:}} When $X$ is a connected Stein manifold of
dimension $n\ge 2$, the above Claim was proved by Laufer (Theorem 9
in \cite{Lauf}).

\subsection{Theorem 1.5 for Stein spaces imply Theorem 1.4 (or Corollary 4.4 in \cite{BS})}

\medskip
\noindent Let $(X, \mathcal{O})$ be a reduced Stein space and $K$ be
a compact set such that $X\setminus K$ has no relatively compact
branches in $X$ and furthermore $\text{codh}(\mathcal{O}_X)\ge 2$.
When $X$ is not irreducible, then we have
$\text{codh}(\mathcal{O}_x)=\text{inf}_{P}\{\text{codh}({\mathcal{O}_{x}})_{
P}\}$, where $m$ is the maximal ideal in $\mathcal{O}_x$,
$\mathcal{P}$ is a prime ideal in $\mathcal{O}_x$ and the infimum is
taken over all prime ideals $\mathcal{P}$ that contain $m$ (see for
example Corollary 3.6, page 42 in \cite{Har}). Let us remark here
that the condition $\text{codh}\,(\mathcal{O}_x)\ge 2$ for all $x\in
K$ (that appears in Corollary 4.4 in \cite{BS}) would guarantee that
$1$-dimensional irreducible components of $X$ do not intersect $K$.

\medskip
\noindent We shall give an alternative proof of Theorem 1.4 (resp.
Corollary 4.4 in \cite{BS}) assuming that the version of Theorem 1.5
for Stein spaces is true. Let $X=\cup X_i$ be the decomposition of
$X$ into irreducible components and let $(\hat{X},\, \xi)$ be a
normalization of $X$. The map $\xi: \hat{X}\to X$ is a one-sheeted
analytic covering map such that
$\hat{\mathcal{O}}_X=\xi_{*}(\mathcal{O}_{\hat{X}})$. Here
$\hat{\mathcal{O}}_X$ is the sheaf of {\it{weakly}} holomorphic
functions on $X$. According to Grauert and Remmert (Global
Decomposition Theorem, page 172 in \cite{GRem}), the irreducible
components of $X$ are the $\xi$-images of the connected components
of $\hat{X}$ and every irreducible component of $X$ contains smooth
points of $X$. We write $\hat{X}=\bigsqcup \hat{X}_i$, i.e.
$\hat{X}$ is the disjoint union of its connected components, which
in turn are the normalizations of the irreducible branches of $X$.
Since $X\setminus K$ has no relatively compact branches in $X$,
$\xi^{-1}(X\setminus K)=\bigsqcup (\hat{X}_i\setminus \hat{X}_i\cap
\xi^{-1}(K))$ will have no relatively compact connected components
in $\hat{X}$. Hence $\hat{X}_i\setminus (\hat{X}_i\cap \xi^{-1}(K))$
has no relatively compact connected components in $\hat{X}_i$, for
every $i$. If we could show that $\hat{X}_i\setminus (\hat{X}_i\cap
\xi^{-1}(K))$ is connected then we could use the version of Theorem
1.5 and finish the proof of Theorem 1.4.

\medskip
\noindent {\it{Some remarks.}} Let $(X, \mathcal{O})$ be any complex
space and $L$ be any compact subset of $X$. Assume that
$X=L\cup\,(\cup_{i}{V_i})\cup \,(\cup_{a} {U_a})$ where $V_i$ are
the relatively compact connected components of $X\setminus L$ and
$U_a$ the unbounded connected components of $X\setminus L$. Any
compact $L'$ that contains $L$ satisfies $(X\setminus L')\cap U_a
\neq \emptyset$ for each $a$. Take as $L'$  the holomorphic
convex-hull of $L$ in $X$. If we could show that $X\setminus L'$ is
connected, then it would easily follow that $X\setminus L\setminus
(\cup_i\,{V_i})$ is connected.

\medskip
\noindent We want to show that $\hat{X}_i\setminus (\hat{X}_i\cap
\xi^{-1}(K))$ is connected; each $\hat{X}_i$ is a connected, Stein,
normal space and those $\hat{X}_i$ that intersect $\xi^{-1}(K)$ are
at least $2$-dimensional. Let us consider the compact subset
$L:=\hat{X}_i\cap \xi^{-1}(K)$ of $\hat{X}_i$. Let $L'$ be the
holomorphic convex-hull of $L$ in $\hat{X}_i$. Based on the above
remarks, to show that $\hat{X}_i\setminus L$ is connected it would
be sufficient to show that $\hat{X}_i\setminus L'$ is connected. If
the latter set had more than one unbounded connected components,
then one could define an $f\in \Gamma(\hat{X}_i\setminus
L',\,\mathcal{O})$ such that $f=1$ on one of the components and $0$
elsewhere on $\hat{X}_i\setminus L'$.  Now, $L'$ is holomorphically
convex in $\hat{X}_i$, $\text{codh}(\mathcal{O}_y)\ge 2$ for all
$y\in L'$ and $\hat{X}_i$ is Stein. By Theorem 2.1 we have that

\begin{equation}\label{eq:conn}
\Gamma(\hat{X}_i, \mathcal{O})\cong \Gamma(\hat{X}_i\setminus L',
\mathcal{O})
\end{equation}

\noindent By $(\ref{eq:conn})$ we conclude that $f$ extends to a
global section $\hat{f}\in \Gamma(\hat{X}_i, \,\mathcal{O})$. Since
$\hat{X}_i$ is connected, Identit\"atssatz  would imply that
$\hat{f}\equiv 0$ on $\hat{X}_i$ which is impossible. Hence we can
conclude that $\hat{X}_i\setminus L'$ is connected and thus
$\hat{X}_i\setminus (\hat{X}_i\cap (\xi^{-1}(K))$ is connected for
all $i$ such that $\hat{X}_i\cap \xi^{-1}(K)\neq \emptyset$.

\medskip
\noindent  All the assumptions of Theorem 1.5 are satisfied for
$\hat{X}_i$ and $\hat{X}_i\cap \xi^{-1}(K)$. Given $s\in
\Gamma(X\setminus K, \mathcal{O}_X)$; its pull-back $s\circ \xi\in
\Gamma(\hat{X}\setminus \xi^{-1}(K), \mathcal{O}_{\hat{X}})$ and its
restriction on $\hat{X}_i\setminus \hat{X}_i\cap \xi^{-1}(K)$ is a
holomorphic section there. According to Theorem 1.5 there exist
$\sigma_i\in \Gamma(\hat{X}_i, \mathcal{O})$ such that
$\sigma_i=s\circ \xi$ on $\hat{X}_i\setminus \hat{X}_i\cap
\xi^{-1}(K)$. Hence, there exists a $\sigma\in \Gamma(\hat{X},
\mathcal{O}_{\hat{X}})$ such that $\sigma=\sigma_i$ on each
$\hat{X}_i$. Let $\tilde{s}=\xi_{*}\sigma$. Clearly, $\tilde{s}$ is
a weakly holomorphic function on $X$ and it is equal to $s$ on
$X\setminus K$. Let $A:=\{x\in X; \;\;\tilde{s}_x\notin
\mathcal{O}_x\;\}=\text{Supp}\left((\mathcal{O}_X+\tilde{s}\,\mathcal{O}_X)/\mathcal{O}_X\right)$.
As $\mathcal{O}_X+\tilde{s}\,\mathcal{O}_X$ is a coherent
$\mathcal{O}_X$-submodule of $\tilde{O}_{X},\;\;A$ is an analytic
set in $X$ which is contained in $K$. Hence $A$ is a finite subset
of $K$. If $A\neq \emptyset$ then $\text{codh}(\mathcal{O}_y)<2$ for
$y\in A$, otherwise $\tilde{s}$ could extend holomorphically in a
neighborhood $U$ of $y$, by Corollary 3.3, page 38 in \cite{BS}. But
this contradicts the assumption that $\text{codh}(\mathcal{O}_x)\ge
2$ for all $x\in K$. Hence $A=\emptyset$, thus $\tilde{s}\in
\Gamma(X, \mathcal{O})$. This proves the surjectivity of the
restriction mapping $\Gamma(X, \mathcal{O})\to \Gamma (X\setminus K,
\mathcal{O})$. To prove the injectivity  we consider $s\in \Gamma(X,
\mathcal{O})$ such that $s=0$ on $X\setminus K$; hence $s=0$ on each
$X_i\setminus (K\cap X_i)$. By the Identit\"atssatz we have that
$s=0$ on each $X_i$, hence $s=0$ on $X$.

\medskip
\noindent {\bf{Proof of Corollary 1.6.}} To prove Corollary 1.6 we
shall use the above approach, passing to the normalization and
assuming Theorem 1.5 for normal Stein spaces.  Under the assumptions
of Corollary 1.6, we can show that the restriction homomorphism
$\Gamma (U,\mathcal{O})\to \Gamma(U\setminus K,\,\mathcal{O})$ is
actually an isomorphism. This will imply $H^i_K(U,\,\mathcal{O})=0$
for $i=0,\,1$. By excision we have $H^i_K (D,\mathcal{O})=0$ for
$i=0,\,1$. Hence, the restriction homomorphism
$\Gamma(D,\,\mathcal{O})\to \Gamma(D\setminus K,\,\mathcal{O})$ is
an isomorphism.

\subsection{Theorem 1.4 implies Theorem 1.5 for Stein spaces}

\medskip
\noindent We assume Theorem 1.4 is true. We will show how to derive
Theorem 1.5 from it in the case of Stein spaces. Let $\Omega$ be any
domain in $X$ and $K\subset \Omega$ be a compact such that
$\Omega\setminus K$ is connected.

\medskip
\noindent{\bf{Claim:}} Under the assumptions of Theorem 1.5 on
$X,\,\Omega, K$ we have that $X\setminus K$ has no relatively
compact connected components and is connected.

\medskip
\noindent{\it{Proof of the Claim:}} Let $\hat{K}$ be the
holomorphically convex hull of $K$ in $X$. We have seen that
$X\setminus \hat{K}$ has no relatively compact connected components
and that $X\setminus \hat{K}$ is connected (by the results in
sections 3.1 and 3.2) consisting only of one unbounded component.
Let $\{V_i\}_{i\in I}$ (resp. $\{U_a\}_{a\in J}$) be the connected
components of $X\setminus K$ that are relatively compact (resp.
unbounded) in $X$. Since $X\setminus \hat{K}$ is connected this will
imply that the set $J$ consists of a single element. Since $\Omega$
is a neighborhood of $K$ it would intersect every component of
$X\setminus K$. We can write $\Omega\setminus K=\Omega\cap
(X\setminus K)=\Omega\cap \left (\cup_{i\in I}\,V_i \cup
U\right)=\left(\cup_{i\in I}\;\Omega\cap V_i\right)\cup (\Omega\cap
U)$. If $I\neq \emptyset$ then we arrive at a contradiction since
$\Omega\setminus K$ is connected. Hence $X\setminus K$ has no
relatively compact connected components and it is connected.

\medskip
\noindent On  a normal space $X$ of dimension $n\ge 2$ we have that
$\text{codh}(\mathcal{O}_X)\ge 2$. Hence all the assumptions of
Theorem 1.4 are satisfied. Let $s\in \Gamma(X\setminus K,
\mathcal{O}_X)$. By theorem 1.4 there exists a unique $\hat{s}\in
\Gamma(X, \mathcal{O}_X)$ such that $\hat{s}=s$ on $X\setminus K$.
This implies that $H^i_K(X,\mathcal{O})=0$ for $i=0,\,1$ and by the
excision property that $H^i_K (\Omega,\,\mathcal{O})=0$ for
$i=0,\,1$. The latter though imply that $H^0(\Omega,
\,\mathcal{O})\cong H^0(\Omega\setminus K,\,\mathcal{O})$.

\subsection{Necessity of conditions in theorems 1.3 or 1.4}

\noindent Let up point out that the homological codimension (depth)
bound in  Theorem 1.3 (or Theorem 1.4) is necessary for a Hartogs
extension theorem to hold. In \cite{Harv}, Harvey gave an example of
an irreducible, Stein {\it{surface}} $X$ in $\mathbb{C}^4$ with a
single isolated singularity at the origin such that when $K=\{0\}$,
there exists a holomorphic function on $X\setminus K$ that can not
be extended holomorphically to $X$. Harvey's point was to show that
simply requiring the singular Stein space to be of dimension bigger
or equal to 2 does not suffice for an extension theorem to hold.
Since this example has been the source for interesting phenomena
(similar examples had been constructed earlier to show that
irreducible Stein varieties whose singular locus has codimension at
least 2 are not necessarily normal (\cite{GunR}, Volume II, page
196), we will recall it and also provide a calculation of the
homological codimension of the structure sheaf of the variety at $0$
in absence of a reference for it.

\medskip
\noindent Let $f: \mathbb{C}^2\to \mathbb{C}^4$ be the map described
by $f(x,y)=(x^2, x^3, y, xy)$. Clearly $f^{-1}(0,0,0,0)=(0,0)$, $f$
is proper and the Jacobian of $f$ has rank $1$ at $0$ and $2$
everywhere else. If we  set $X:=f(\mathbb{C}^2)$ then  $X$ is a
Stein surface in $\mathbb{C}^4$ with an isolated singularity at $0$
and since $X\setminus \{0\}=f(\mathbb{C}^2\setminus (0,0))$ is
connected $X$ is irreducible. Let us consider the following function
on $X\setminus \{0\}$; $g(z_1,z_2,z_3,z_4)=\frac{z_2}{z_1}=x$, if
$z_1\neq 0$ and $g(z_1,z_2,z_3,z_4)=\frac{z_4}{z_3}=x$, if $z_3\neq
0$. Clearly $g$ is a holomorphic function on $X\setminus\{0\}$. If
$g$ could be extended holomorphically across $0$ to some function
$\tilde{g}$, then one could define a holomorphic map $F: X\to
\mathbb{C}^2$ such that
$F(z_1,z_2,z_3,z_4)=(\tilde{g}(z_1,z_2,z_3,z_4), z_3)$ that is the
inverse of $f$, i.e. $f$ would be a biholomorphism. Hence $g$ can
not be extended holomorphically across $0$ and therefore
$\text{codh}(\mathcal{O}_{_{X,\,0}})\le 1$. We shall show that
$\text{codh}(\mathcal{O}_{_{X,\,0}})=1$ by showing that
$\text{codh}(\mathcal{O}_{_{X,\,0}})>0$. We know that
$\text{codh}(\mathcal{O}_{_{X,\,0}})=0$ if and only if every element
in $m_{_{X,\,0}}$ (the maximal ideal in $\mathcal{O}_{_{X,\,0}})$ is
a zero-divisor in $\mathcal{O}_{_{X,\,0}}$. Taking into account that
$\mathcal{O}_{_{X,\,0}}\cong \mathbb{C}\{x^2,\,x^3,\,y,\,xy\}$ and
that $m_{_{X,\,0}}$ is generated by elements in the maximal ideal of
$\mathcal{O}_{\mathbb{C}^4}$ at $0$, we immediately see that
$\text{codh}(\mathcal{O}_{_{X,\,0}})$ can not be zero.

\medskip
\noindent Let $\nu$ be the minimal depth (homological codimension)
of the rings $\mathcal{O}_{X,x}$, $x\in \text{Sing}X$. When $\nu=1$
Andersson and Samuelsson gave a necessary and sufficient condition
for extendability. More precisely they proved the following:

\begin{theorem}(Theorem 1.4 ii) in \cite{AS}) Assume that $X$ is a
Stein space of pure dimension $d$ with globally irreducible
components $X^{\ell}$ and let $K$ be a compact subset such that
$\text{Reg}\,X^{\ell}\setminus K$ is connected for each $\ell$.
Assume $\nu=1$ and let $\chi$ be a cut-off function that is
identically equal to $1$ in a neighborhood of $K$. There exists a
smooth $(d,\,d-1)$ form $a$ on $\text{Reg}\,X$ such that the
function $\phi\in \mathcal{O}(X\setminus K)$ has a holomorphic
extension $\Phi$ across $K$ if and only if

\begin{equation}\label{eq:mc}
\int_{X} \overline\partial \chi\wedge a\,\phi\,h=0
\end{equation}

\noindent for all $h\in \mathcal{O}(X),$ and where the integrals
exist as principal values at $\text{Sing}\,X$.
\end{theorem}

\medskip
\noindent In the case of curves in $\mathbb{C}^2$ similar moment
conditions were given by Hatziafratis (Theorem 2 in \cite{Hatz}).

\medskip
\noindent {\bf{Remark:}} Using Harvey's example we can show why a
stronger condition on the homological codimension of $\mathcal{O}_X$
is needed in the statement of Corollary 4.4 in \cite{BS}. There, it
is only assumed that $\text{codh}(\mathcal{O}_x)\ge 2$ for all $x\in
K$. In the proof provided though, one has to use the fact that
$\text{codh}(\mathcal{O}_x)\ge 2$ for all $x\in \hat{K}_X$, the
holomorphic convex hull of $K$ in $X$. We shall show that we can
find a compact subset $K$ in Harvey's surface $X$ such that $X
\setminus K$ is connected and $0\in \hat{K}_X$ (recall that
$\text{codh}(\mathcal{O}_0)=1$!). Let $P(r)$ denote the bi-disc of
radius $r$. We can construct a compact $L$ in $\mathbb{C}^2$ such
that: $L:=\overline{P(1)}\setminus T_{\epsilon}$, where
$T_{\epsilon}:=\{(z,w);\;\;(\text{Im}\,z)^2+|w|^2<\epsilon^2\}$ for
$0<\epsilon<1$ is an open tube about $\mathbb{R}\times\{0\}$.
Clearly $0\notin L$ and $\mathbb{C}^2\setminus L$ is connected. In
what follows by $\Delta$ we will denote the unit disk in
$\mathbb{C}$ with respect to the $z$ or $w$ variable.

\smallskip
\noindent{\bf{Claim:}} If $g$ is a holomorphic function near $L$
then $g_{\upharpoonright L}$ has a bounded holomorphic extension
$\tilde{g}$ to $P(1)$.

\smallskip
\noindent{\it{Proof of Claim.}} Let us set
$\tilde{g}(z,w):=\frac{1}{2\pi\, i}\;\int_{|\zeta|=1}
\dfrac{g(z,\zeta)}{\zeta-w}\,d\zeta$. Clearly $\tilde{g}$ is defined
in $P(1)$. When $|z|<1$ and $|\text{Im}\,z|>\epsilon$, we have
$\{z\}\times \overline{\Delta}\subset L$, so $\tilde{g}(z,w)=g(z,w)$
when $|w|<1$ by Cauchy's formula. Since $L$ is connected and
$U:=(\{|z|<1,\;|\text{Im}\,z|> \epsilon\})\times \{|w|<1\}\neq
\emptyset$ subset of $L$,  by the identity principle of holomorphic
functions we conclude that $\tilde{g}=g$ on $L$ and hence
$\tilde{g}$ is a bounded holomorphic extension of $g$ in $P(1)$.

\medskip
\noindent Let $f:\mathbb{C}^2\to X$ be the parametrization of
Harvey's example and let $K:=f(L)$. If $g$ is a holomorphic function
on a neighborhood $D$ of $K$ in $X$, it follows that $g$ has a
bounded analytic continuation on $f(P(1))\setminus \{0\}$ since $f$
is a biholomorphism from $\mathbb{C}^2\setminus \{(0,0)\}$ to
$X\setminus \{0\}$. Hence, $D$ can not be a Stein domain in
$X\setminus \{0\}$ and $K$ can not be contained in a holomorphically
convex compact in $X\setminus \{0\}$.

\subsection{Relatively compact branches} Take $X=:\mathbb{C}^2$ and
$K=:\{(z,w)\;|\;1\le |z|^2+|w|^2\le 2\}$. Then
$\text{codh}(\mathcal{O}_X)=2=\text{dim}\,X$, $X$ is normal
connected (thus irreducible) and $X\setminus K=U\cup W$ where $U$ is
the unit ball and $W=\{(z,w)\;|\;|z|^2+|w|^2>2\}$. Consider a
function $f$ such that  $f=1$ on $U$ and $f=0$ on $W$. Clearly $f\in
H^0(X\setminus K, \mathcal{O})$ but $f$ can not be extended to a
holomorphic function on $\mathbb{C}^2$. One can produce more refined
counterexamples where $X\setminus K$ is connected; consider for
example as $X:=X_1\cup X_2$ where $X_1$ is the plane in
$\mathbb{C}^3$ given by $z_1=0$ and $X_2$ the plane in
$\mathbb{C}^3$ described by $z_3=0$ and let
$K:=\{(z_1,z_2,z_3);\;1\le |z_1|^2+|z_2|^2\le 2, \,z_3=0\;\}$.
Clearly $\text{codh}\,(\mathcal{O}_X)=2$ and $X\setminus K$ is
connected. On the other hand $X\setminus K$ has a relatively compact
component $U$. Let us define a function $f$ such that $f=z_1$ on $U$
and $f=0$ elsewhere on $X\setminus K$; $f\in \Gamma(X\setminus
K,\;\mathcal{O})$. Such an $f$ can not be extended holomorphically
on $X$. In \cite{AS} an example is given (Example 3, Section 8) to
show that when $X\setminus K$ has an irreducible component that is
relatively compact in $X$ and $\nu=1$ the moment condition
$(\ref{eq:mc})$ is not sufficient to guarantee extendability.

\section{A new proof of a Hartogs extension theorem based on Andreotti-Grauert's work}

\medskip
\noindent In \cite{AG} the following theorem was proved:

\begin{theorem} (Th$\acute e$or$\grave e$me 15) Let $X$ be a complex
space for which there exists a function $\phi>0$, strongly
$q$-convex such that the subsets $X_{\epsilon,\,c}:=\{x\in
X\;|\;\;\epsilon<\phi(x)<c\}$ are relatively compact for
$\epsilon,\,c>0$. Let $B_c:=\{x\in X\;|\;\;\phi(x)<c\}$ pour $c>0$.
Let $\mathcal{F}$ be a coherent sheaf on $X$. Then the homomorphism

$$
H^r (X,\,\mathcal{F})\to H^r(X\setminus B_c,\,\mathcal{F})
$$

\noindent is bijective for $r<\text{codh}(\mathcal{F})-q$ and
injective for $r=\text{codh}(\mathcal{F})-q$.
\end{theorem}

\medskip
\noindent The key ingredient in the proof is the vanishing of the
cohomology groups with compact support, $H^r_c (B_c, \mathcal{F})$
for all $r\le \text{codh}(\mathcal{F})-q$. This latter vanishing is
obtained via Grauert's bumbing method and some local vanishing of
the sheaf cohomology groups with compact support near the boundary
of $B_c$. When $X$ is a Stein space (or equivalently a $1$-complete
space) such that $\text{codh}(\mathcal{O}_X)\ge 2$ then from the
above theorem we deduce that the restriction homomorphism

$$
H^0(X,\,\mathcal{O})\to H^0(X\setminus B_c,\,\mathcal{O})
$$

\noindent is a bijection. From this we can deduce the following
theorems

\begin{theorem} Let $X$ be a Stein space such that
$\text{codh}(\mathcal{O}_X)\ge 2$ and $K$ be a compact subset of
$X$. If $X\setminus K$ has no irreducible components that are
relatively compact in $X$ then the restriction homomorphism
$H^0(X,\,\mathcal{O})\to H^0(X\setminus K,\,\mathcal{O})$ is an
isomorphism.
\end{theorem}

\noindent{\it{Proof.}} Since $X$ is Stein, there exists a strictly
plurisubharmonic function $\phi$ that satisfies the assumptions of
Theorem 15 in \cite{AG} and for $K$ compact subset of $X$ there
exists a $c>0$ such that $K\subset B_c$. Let $f\in H^0(X\setminus
K,\,\mathcal{O})$. Then the restriction of $f$ on  $X\setminus B_c$
is a holomorphic section there, and by Theorem 4.1 there exists a
unique section $\hat{f}\in H^0(X,\mathcal{O})$ such that $\hat{f}=f$
near  $X\setminus B_c$. Now by assumption each irreducible component
$U$ of $X\setminus K$ is unbounded, hence it will intersect
$X\setminus {\overline{B}_c}$\; i.e.\; $U\cap
X\setminus{\overline{B}_c}\neq \emptyset$. We have $\hat{f}=f$ on
the non-empty open subset $U\cap (X\setminus \overline{B}_c)$ of
$U$. By the identity principle for holomorphic functions  we obtain
$\hat{f}=f$ on $U$ and hence $\hat{f}=f$ on $X\setminus K$. As a
consequence we obtain that $H^i_K(X,\,\mathcal{O})=0$ for $i=0,\,1$.

\noindent \begin{theorem} Let $X,K$ be as above and let $D$ be a
neighborhood of $K$ in $X$. Then the restriction homomorphism
$H^0(D,\,\mathcal{O})\to H^0(D\setminus K,\,\mathcal{O})$ is an
isomorphism.
\end{theorem}

\noindent {\it{Proof.}} As a consequence of the previous theorem we
obtained that $H^i_K(X,\,\mathcal{O})=0$ for $i=0,\,1$.  By excision
we also have that $H^i_K(D,\,\mathcal{O})=0$ for $i=0,\,1$. The
vanishing of $H^1_K(D,\,\mathcal{O})$ implies the surjectivity of
$H^0(D,\,\mathcal{O})\to H^0(D\setminus K,\,\mathcal{O})$. The
injectivity of the restriction homomorphism follows from the
exactness of $0\to H^0_K(D,\,\mathcal{O})\to H^0(D,\,\mathcal{O})\to
H^0(D\setminus K,\,\mathcal{O})\to H^1_K(D,\,\mathcal{O})=0$ and the
fact that $H^0_K(D,\mathcal{O})=0$.

\medskip
\noindent {\bf{Remark:}} As we have seen in section 3, Theorem 4.3
implies Theorem 1.5 for Stein spaces which in its turn implies
Corollary 4.4 in \cite{BS}.

\medskip
\noindent Using these theorems we can easily prove connectedness of
the complements of some compact sets in appropriate Stein spaces.
More precisely we obtain the following:

\begin{corollary} Let $X$ be an irreducible Stein space with
$\text{codh}(\mathcal{O}_X)\ge 2$. Let $K$ be a compact subset of
$X$ such that $X\setminus K$ has no irreducible components that are
relatively compact in $X$. Then $X\setminus K$ is connected.
\end{corollary}

\noindent {\it{Proof.}} Suppose on the contrary that $X\setminus
K=U\cup W$ where $U,\,W$ disjoint open sets. Let $f=1$ on $U$ and
$f=0$ on $W$. Clearly $f\in H^0(X\setminus K,\,\mathcal{O})$. By
Theorem 4.2 there exists $\hat{f}\in H^0(X,\,\mathcal{O})$ such that
$\hat{f}=f$ on $X\setminus K$. By the identity principle we could
have that $\hat{f}=1$ on $X$ which is not possible.

\begin{corollary} Let $X,\,K$ be as in Theorem 4.3 and let $D$ be a connected
neighborhood of $K$. Then $D\setminus K$ is connected.
\end{corollary}

\noindent{\it{Proof.}} Follows along the same lines as the proof of
Corollary 4.4 using Theorem 4.3 and the identity principle.

\medskip
\noindent{\bf{Remark:}} If $D$ is a neighborhood of $K$, then $D$
intersects every component of $X\setminus K$; hence $D\setminus K$
connected would imply $X\setminus K$ connected.

\section{A $\overline\partial$-approach to Hartogs extension theorems on Stein spaces}

\medskip
\noindent The main goal in this section is to streamline an analytic
approach via $\overline\partial$ techniques and duality theorems for
proving Theorem 1.5 for Stein spaces with arbitrary singularities.

\medskip
\noindent First we need some preliminaries on the existence of
smooth cut-off functions on a singular complex space.

\subsection{Cut-off functions.} A function on a complex space $X$ is
$C^{\infty}$ iff it is locally the pull-back of a
$C^{\infty}$-function by a local embedding in an open set in some
$\mathbb{C}^N$. This condition is independent on the choice of the
embedding. In the case where there is a proper embedding $\phi: W\to
U^{\text{open}}\subset \mathbb{C}^N$ we may choose a cut-off
function $\chi'\in C^{\infty}_0(U)$ with $\chi'=1$ on
$\overline{\phi(V)}$ for $V\Subset W$ and set $\chi:=\chi'\circ
\phi$.

\medskip
\noindent When $X$ is a complex space countable at infinity but not
necessarily Stein, the proof of the existence of a smooth partition
of unity subordinate to an open covering carries over from the
manifold case. This implies the existence of smooth cut-off
functions in general.

\subsection{The analytic approach; a reduction.}

\medskip
\noindent Let $(X, \mathcal{O})$ be a connected, normal, Stein space
of dimension greater or equal to $2$. Let $K$ be a compact subset of
$X$ and $D$ an  open neighborhood of $K$ in $X$ such that
$D\setminus K$ is connected. We need to show that given any $s\in
\Gamma(D\setminus K, \mathcal{O})$ there exists unique $\hat{s}\in
\Gamma(D, \mathcal{O})$ such that $\hat{s}_{|_{D\setminus K}}=s$. We
shall need the following Proposition:

\begin{proposition} Let $X$ be a connected, non-compact normal
complex space and $K$ a compact subset of $X$. Let us assume that
$K$ has an open neighborhood $\Omega$ in $X$ with the following
property (P)

\medskip
\noindent For every ``nice'' $\overline\partial$-closed $(0,1)$-form
$f$ on $\text{Reg}\,\Omega$ with
$\text{supp}_X\,f:=\text{closure}_X\,(\text{supp}\,f)$ compact in
$\Omega$, the equation $\overline\partial u=f$ has a solution on
$\text{Reg}\,\Omega$ with
$\text{supp}_X\,u:=\text{closure}_X\,(\text{supp}\,u)$ compact in
$\Omega$.

\medskip
\noindent Then, for every open  neighborhood $D$ of $K$ with
$D\setminus K$ connected and every $s\in \Gamma(D\setminus
K,\,\mathcal{O})$ there exists a unique $\tilde{s}\in
\Gamma(D,\,\mathcal{O})$ such that $\tilde{s}=s$ on $D\setminus K$.
\end{proposition}

\medskip
\noindent ``Nice'' in the above statement means that the form $f$
can be smoothly extended on $U$ for some local embedding $\phi:
W\subset \Omega \to U^{\text{open}}\subset \mathbb{C}^N$.

\medskip
\noindent{\it{Proof:}} We choose $\chi\in C^{\infty}_o (D\cap
\Omega)$ such that $\chi=1$ near $K$ and we define a $(0,1)$ form
$f$ on $\text{Reg}\,\Omega$ as follows:

$$
f=\left\{\begin{array}{ccc} s\,\overline\partial \chi &\;
\text{on}\;& \text{supp}(\overline\partial \chi)\cap
\text{Reg}\,\Omega,
\\ 0 &\;\text{on}\;& \text{Reg}\Omega\setminus
\text{supp}(\overline\partial \chi)
\\
\end{array}\right.
$$

\medskip
\noindent Let  $u$ be a solution to $\overline\partial u=f$ on
$\text{Reg}\,\Omega$ with $\text{supp}_X u$ compact in $\Omega$. Let
$u^0,\;f^0$ denote the trivial extensions by zero outside
$X\setminus \Omega$ of $u$ and $f$ respectively. Then
$\overline\partial u^0=f^0$ on $\text{Reg}\,X$ and
$\overline\partial u^0=0$ on $\text{Reg}\,X\setminus\text{supp}\,f$.
Let $U$ be an unbounded connected component of $X\setminus
\text{supp}\,\chi$. Clearly $u^0$ is holomorphic on
$\text{Reg}\,X\setminus \text{supp}\chi$, so $u^0$ must be
identically equal to $0$ on $U$. Hence $\text{supp}\,u^0\subset
X\setminus U\cap \Omega$.

\medskip
\noindent Let us define $\tilde{s}$ on $\text{Reg}\,D$ by setting
$\tilde{s}:=u^0+(1-\chi)\,s$. Clearly $\overline\partial
\tilde{s}=0$ on $\text{Reg}\,D$ and  extends as a holomorphic
section (still denoted by $\tilde{s}$) on $D$ by normality. Since
$X$ is connected, $U$ must have nonempty boundary $bU$  which must
intersect the boundary of $\text{supp}\,\chi$. Hence $U\cap
(D\setminus K)\neq \emptyset$ and thus $\tilde{s}=s$ on this open
set. Since $D\setminus K$ is connected normal (thus irreducible) we
obtain $\tilde{s}=s$ on $D\setminus K$ by the identity principle.
Uniqueness follows from the fact that $D$ is connected and the
identity principle for holomorphic functions.

\medskip
\noindent Hence it remains to prove that we can always find an open
neighborhood $\Omega$ of $K$ that satisfies property (P) of the
above proposition. Since $K$ is compact on a Stein space we can
always find an open, relatively compact Stein subdomain of $X$ that
contains $K$. We shall establish below that Property (P) holds for
such domains. 

\subsection{Solvability with compact support.}
\medskip
\noindent Let us recall the settings and the main theorem in
\cite{FOV1}. Let $X$ be a pure $n-$dimensional reduced Stein space,
$A$ a lower dimensional complex analytic subset with empty interior
containing $\text{Sing}\,X$, the singular locus of $X$. Let $\Omega$
be an open relatively compact Stein domain in $X$ and
$L=\widehat{\overline{\Omega}_X}$ be the holomorphic convex hull of
$\overline{\Omega}$ in $X$.  Since $X$ is Stein and
$L=\widehat{L_X},$  $L$ has a neighborhood basis of Oka-Weil domains
in $X$. Let $X_0$ be an Oka-Weil neighborhood of $L$ in $X,\;
X_0\subset\subset X$. Then $X_0$ can be realized as a holomorphic
subvariety of an open polydisk in some $\mathbb{C}^s$. Set
$\Omega^*=\Omega\setminus A$. Let $d_A$ be the distance to $A$,
relative to an embedding of $X_0$ in $\mathbb C^s$ and let $|\;|$
and $dV$ denote the induced norm on $\Lambda^{\cdot} T^*\Omega^*,$
resp. the volume element (different embeddings of neighborhoods of
$\overline{\Omega}$ in $\mathbb{C}^s$ give rise to equivalent
distance functions and norms).  For a  measurable $(p,q)$ form $u$
on $\Omega^*$  let $\|u\|^2_{N,\Omega}:=\int_{\Omega^*}|u|^2
d_A^{-N} dV$.

\medskip
\noindent In \cite{FOV1} we proved the following theorem:

\begin{theorem}
Let $X,\,\Omega$ be as above. For every $N_0 \geq 0,$ there exists
$N \geq 0$ such that if $f$ is a
 $\overline{\partial}-$closed $(p,q)-$form on $\Omega^*,q>0,$
 with $\|f\|_{N,\Omega}<\infty,$ there is  $v \in L^{2,\,{\rm loc}}_{p,q-1}(\Omega^*)$ solving $\overline{\partial}v=f,$ with $\|v\|_{N_0,\Omega'}<\infty$
 for every $\Omega' \subset \subset \Omega$. For each $\Omega'\subset
 \subset \Omega,$ there is a solution of this kind satisfying
 $\|v\|_{N_0,\Omega'}\leq C \|f\|_{N,\Omega},$ where $C$ is a positive
constant that depends only on
 $\Omega', N,N_0.$
 \end{theorem}

\medskip
\noindent  We shall establish a dual version of Theorem 5.2 taking
as $A:=\text{Sing}\,X$. We would like to prove the following theorem

\begin{theorem} Let $f$ be a
$(p,q)$ form defined on $\text{Reg}\,\Omega$ and
$\overline\partial$-closed there with $0<q<n$, compactly supported
in $\Omega$ and such that $\int_{\text{Reg}\,\Omega} |f|^2\,
d_A^{N_0} dV<\infty$ for some $N_0\ge 0$. Then there exists a
solution $u$ to $\overline\partial u=f$ on $\text{Reg}\,\Omega$
satisfying $\text{supp}_X\, u\Subset \Omega$ and such that

$$\int_{\text{Reg}\,\Omega} |u|^2\, d_A^{N} dV\le
C\,\int_{\text{Reg}\,\Omega} |f|^2\, d_A^{N_0} dV
$$

\noindent where $N$ is a positive integer that depends on $N_0$ and
$\Omega$ and $C$ is a positive constant that depends on
$N_0,\,N,\,\Omega$\; and\;\;$\text{supp}\,f$.

\end{theorem}

\medskip
\noindent {\it{Proof of theorem 5.3.}} The proof is inspired by the
proof of Serre's duality theorem on complex manifolds. We shall
apply Theorem 5.2 twice for forms of different bidegree and for
$A=\text{Sing}X$. Let $N_0$ be as in Theorem 5.3. Then according to
Theorem 5.2 for the pair $(X,\,\Omega)$ and for $N_0+2$ there exists
$N_1 (>>N_0+2)$ such that if $F\in L^{2,\,loc}_{n-p,\,
n-q}(\text{Reg} \,\Omega)$, $\overline\partial F=0$ on
$\text{Reg}\,\Omega$ with $\|F\|_{N_1,\, \Omega}<\infty$ then, there
exists $a\in L^{2,\,loc}_{n-p,\,n-q-1}(\text{Reg}\,\Omega)$
satisfying $\overline\partial a=F$ on $\text{Reg}\,\Omega$ and such
that $\|a\|_{N_0+2,\, Y}<\infty$ for all $Y \Subset \Omega$.
Similarly, for the above $N_1$ there exists an $N (>> N_1)$ such
that if $G\in
L^{2,\,loc}_{n-p,\,n-q+1}(\text{Reg}\,\Omega),\;\overline\partial$-closed
on $\text{Reg}\,\Omega$ with $\|G\|_{N,\,\Omega}<\infty$ then, there
exists a \;\, $b\in L^{2,\,loc}_{n-p,\,n-q}(\text{Reg}\,\Omega)$
satisfying $\overline\partial b=G$ on $\text{Reg}\,\Omega$ and such
that $\|b\|_{N_1,\,Y}<\infty$ for all $Y\Subset \Omega$.

\medskip
\noindent Let $f,\,N_0$ be as in Theorem 5.3 and $N,\,N_1$ be chosen
as above.  Let us define

$$\mathcal{H}_N:=\{w\in L^{2,loc}_{n-p,\,n-q+1}(\text{Reg}\,\Omega):\;\;
\|w\|_{N,\, \Omega}<\infty\}$$

\medskip
\noindent and $\mathcal{H}^{loc}_{N_1}:=\{v\in
L^{2,\,loc}_{n-p,\,n-q}(\text{Reg}\,\Omega):
\;\|v\|_{N_1,\,Y}<\infty\; \; \text{for\;\;all\;\;}\;\; Y\Subset
\Omega\}$. Consider the following map:

$$ L_f: \mathcal{H}_N\cap \text{kern}(\overline\partial)\to
\mathbb{C}$$

\noindent defined by $$L_f(w):=(-1)^{p+q+1}\;\int_{Reg\,\Omega}
v\wedge f$$

\noindent where $v\in \mathcal{H}^{loc}_{N_1}$ is a solution to
$\overline\partial v=w$ on $\text{Reg}\,\Omega$ (such a solution
always exists by Theorem 5.2).

\medskip
\noindent First of all we need to show that $L_f$ is well-defined,
i.e. independent of the choice of the solution $v\in
\mathcal{H}^{\text{loc}}_{N_1}$ to the equation $\overline\partial
v=w$ on $\text{Reg} \Omega$. It suffices to show that $\int_{Reg
\Omega} v\wedge f=0$ when  $v\in \mathcal{H}^{loc}_{N_1}$ and
$\overline\partial v=0$ on $\text{Reg}\,\Omega$. Clearly
$\overline\partial v=0$ on $\text{Reg}\,\Omega$ and
$\|v\|_{N_1,\,V}<\infty$. Hence according to what was discussed
above there exists an $a\in
L^{2,\,loc}_{n-p,\,n-q-1}(\text{Reg}\,\Omega)$ satisfying
$\overline\partial a=v$ on $\text{Reg}\,\Omega$ and such that
$\|a\|_{N_0+2,\, Y}<\infty$ for all $Y\Subset \Omega$.

\medskip
\noindent Let $\chi_{\delta}\in C^{\infty}(\text{Reg}\,\Omega)$ such
that $0\le \chi_{\delta}\le 1,$\; $\chi_{\delta}(z)=1$ when
$d_A(z)>\delta,\;\;\;\chi_{\delta}(z)=0$ if $d_A(z)\le
\frac{\delta}{2}$ and such that $|\overline\partial
\chi_{\delta}|\le \dfrac{C}{\delta}$ on $\text{Reg}\,\Omega$,  for
some positive constant $C$ independent of $\delta$. Then
$L_f(z)=\underset{\delta\to 0}\lim
\int_{_{\text{Reg}\,\Omega}}\chi_{\delta}\, v\wedge
f=\underset{\delta\to 0}\lim \int_{_{\text{Reg}\,\Omega}}
\chi_{\delta}\,\overline\partial a\,\wedge f=-\underset{\delta\to
0}\lim \int_{_{\text{Reg}\,\Omega}} \overline\partial
\chi_{\delta}\wedge a\wedge f$. The last equality follows from
integration by parts arguments; the form $\chi_{\delta}\,f$ has
compact support on $\text{Reg}\,\Omega$ but is not really smooth, so
a smoothing argument is needed in order to apply the standard
Stokes' theorem. Let $I_{\delta}:= \int_{_{\text{Reg}\,\Omega}}
\overline\partial \chi_{\delta}\wedge a\wedge f$. By Cauchy-Schwarz
inequality we have

\begin{equation*}
|I_{\delta}|^2\le \left(\int_{\text{supp}\,f\cap
\text{supp}\,(\overline\partial\chi_{\delta})}\;
|d_A\,\overline\partial \chi_{\delta}|^2\;|f|^2\;
d_A^{N_0}\;dV\right)\;\left(\int_{\{{\chi_{\delta}}<1\}\cap
\text{supp}\,f}\; |a|^2\;d_A^{-(N_0+2)}\;dV\right)=A\cdot B
\end{equation*}

\smallskip
\noindent But $B$ is finite and  $A\to 0$ as $\delta\to 0$ since
$|d_A\,\overline\partial \chi_{\delta}|\le C$ and
$\int_{\text{Reg}\,\Omega}\;|f|^2\, d_A^{N_0}\;dV<\infty$. Hence
$\underset{\delta\to 0}\lim I_{\delta}=0$ and thus $L_f$ is
well-defined.

\medskip
\noindent Clearly $L_f$ is a linear map. By Cauchy-Schwarz we have

\begin{equation}\label{eq:CS1}
|L_f(w)|\le C_0\,\left(\int_{\text{supp}\,f} |v|^2\,
d_A^{-N_1}\;dV\right)^{\frac{1}{2}}\;\left(\int_{\text{Reg}\,\Omega}\,|f|^2
\;d_A^{N_0}\;dV\right)^{\frac{1}{2}}
\end{equation}

\noindent To obtain this inequality we used the fact that
$N_1>>N_0$ and hence $d_A^{N_1-N_0}\le C_0$ on
$\overline{\Omega}$. We want to show that $L_f$ factors into a
bounded linear functional on a subspace $\mathcal{A}$ of $H_N\cap
\text{kern}(\overline\partial)$. Recall that  the following lemma
was proven in \cite{FOV1} using the open mapping theorem for
Fr\'echet spaces.

\begin{lemma}(Lemma 4.2 in \cite{FOV1}) Let $M$ be a complex manifold and let $E$ and $F$ be Fr\'echet
spaces of differential forms (or currents) of type $(p,q-1),
\;(p,q)$, whose topologies are finer (possibly strictly finer) than
the weak topology of currents. Assume that for every $f\in F$, the
equation $\overline\partial u=f$ has a solution $u\in E$. Then, for
every continuous seminorm $p$ on $E$, there is a continuous seminorm
$q$ on $F$ such that the equation $\overline\partial u=f$ has a
solution with $p(u)\le q(f)$ for every $f\in F,\;q(f)>0$.
\end{lemma}

\medskip
\noindent Let $p(v):=\left(\int_{\text{supp}\,f} |v|^2\,
d_A^{-N_1}\;dV\right)^{\frac{1}{2}}$. Using the lemma in our
situation, given the seminorm $p$ there exists an open set
$W\Subset \Omega$ (that depends on the $\text{supp}\,f$) such that
for all $w\in \mathcal{H}_N\cap \text{kern}(\overline\partial)$
with $q(w)=\left(\int_W
|w|^2\,d_A^{-N}\,dV\right)^{\frac{1}{2}}>0$ there exists a
solution $v$ to $\overline\partial v=w$ on $\text{Reg}\,\Omega$
and a positive constant $C$ satisfying

\begin{equation}\label{eq:CS2}
\left(\int_{\text{supp}\,f} |v|^2\,
d_A^{-N_1}\;dV\right)^{\frac{1}{2}}\le C \,q(w).
\end{equation}

\medskip
\noindent If $q(w)=0$ then the same argument will imply that for
every $\epsilon>0$ there exists a solution $v_{\epsilon}$ to
$\overline\partial v_{\epsilon}=w$ on $\text{Reg}\,\Omega$ with
$p(v_{\epsilon})<\epsilon$. Then for such a $w$ we have:
$|L_f(w)|\le C_0 \epsilon \int_{\text{Reg}\,\Omega}
|f|^2\,d_A^{N_0}\,dV$. Here we used the fact that $L_f(w)$ is
well-defined independent of the choice of the solution $v\in
\mathcal{H}^{n-q,\,loc}_{N_1}(\Omega)$ to $\overline\partial v=w$.
Taking the limit as $\epsilon\to 0$ we get that

\begin{equation}\label{eq:CS3}
|L_f(w)|=0\le C_0\,q(w)\,\left(\int_{\text{Reg}\,\Omega}\,|f|^2
\;d_A^{N_0}\;dV\right)^{\frac{1}{2}}.
\end{equation}

\medskip
\noindent Combining
$(\ref{eq:CS1}),\,(\ref{eq:CS2}),\,(\ref{eq:CS3})$ we obtain for all
$w\in \mathcal{H}_N\cap \text{kern}(\overline\partial)$

\begin{equation}\label{eq:cont}
|L_f(w)|\le C_0\,C\,q(w)\,\left(\int_{\text{Reg}\,\Omega}\,|f|^2
\;d_A^{N_0}\;dV\right)^{\frac{1}{2}}.
\end{equation}

\medskip
\noindent From $(\ref{eq:cont})$ we see that  $L_f(w)$ depends only
on $w_{\upharpoonright W}$. Indeed, let $w,\,w'\in \mathcal{H}_N\cap
\text{kern}(\overline\partial)$ such that $w_{\upharpoonright
W}=w'_{\upharpoonright W}$. Then
$L_f(w)=L_f(w-w'+w')=L_f(w-w')+L_f(w')$. From $(\ref{eq:cont})$ we
obtain that $|L_f(w-w')|\le C_0\,C\,
q(w-w')\,\left(\int_{\text{Reg}\,\Omega}\,|f|^2
\;d_A^{N_0}\;dV\right)^{\frac{1}{2}}$. But $q(w-w')=0$ as $w-w'=0$
on $W$. Hence $L_f(w-w')=0$ and thus $L_f(w)=L_f(w')$. Hence $L_f$
factors to a well-defined bounded linear functional on
$\mathcal{A}:=\{w_{\upharpoonright W};\;\;w\in \mathcal{H}_N\cap
\text{kern}(\overline\partial)\}\subset \mathcal{H}_N(W)$. Here
$\mathcal{H}_N(W):=\{w\in L^{2,\,loc}_{n-p,\,n-q+1}(W):
\;\|w\|_{N,W}<\infty\}$.

\medskip
\noindent We make a norm-preserving extension of the above
functional $L_f$ to $\mathcal{H}_N(W)$. Let us call $\tilde{L}_f$
the extended functional. By Riesz representation theorem there
exists a $u'\in \mathcal{H}_N(W)$ such that for all $w\in
\mathcal{H}_N(W)$ we have

\begin{equation}\label{eq:Dfn}
\tilde{L}_f(w)=\int_{\text{Reg}\,W} <w,\,u'>\;d_A^{-N}\;dV.
\end{equation}

\noindent Set $u:=d_A^{-N}\;\overline{*}\;u'$ on $\text{Reg}\,W$
(here $*$ is the Hodge-star operator) and extend by zero outside
$\overline{W}$. We claim that $u$ is the desired solution of Theorem
5.3. Certainly $\text{supp}\,u\Subset \Omega$ and
$\int_{\text{Reg}\,\Omega} |u|^2\,d_A^N\;dV=\int_{\text{Reg}\,W}
|u'|^2\; d_A^{-N}\, dV<\infty$, since $u'\in \mathcal{H}_N(W)$. We
can control the weighted $L^2$ norm of $u'$ in terms of the weighted
$L^2$-norm of $f$ taking into consideration the following:

\begin{equation}\label{eq:L2normu'}
\left(\int_{\text{Reg}\,W}
|u'|^2\,d_A^{-N}\,dV\right)^{\frac{1}{2}}= \|\tilde{L}_f\|=
\|{L_{f}}_{\upharpoonright
\mathcal{A}}\|=\phantom{sdsdfsdsdsdfhjjlklliuiu}
\end{equation}

\begin{equation*}
=\text{sup}\{ |L_f(w)|\;: w\in \mathcal{A} \;\text{with}\; q(w)\le 1
\}\le C_0\,C\,\left(\int_{\text{Reg}\,\Omega}\,|f|^2
\;d_A^{N_0}\;dV\right)^{\frac{1}{2}}.
\end{equation*}

\medskip
\noindent
 It remains to show that
$\overline\partial u=f$ on $\text{Reg}\,\Omega$. Let $\phi\in
C^{\infty}_{0,\,(n-p,\, n-q)}(\text{Reg}\,\Omega)$ be a smooth
compactly supported form of bidegree $(n-p,\,n-q)$ on
$\text{Reg}\,\Omega$. We need to show that

\begin{equation}\label{eq:Lsolv}
\int_{\text{Reg}\,\Omega} \overline\partial \phi\wedge u
=(-1)^{p+q+1}\, \int_{\text{Reg}\,\Omega}  \phi \wedge f.
\end{equation}

\medskip
\noindent But $\phi \in \mathcal{H}^{loc}_{N_1}$,
$\overline\partial\phi \in \mathcal{H}_N$ and $\overline\partial
\phi_{\upharpoonright W}\in \mathcal{A}$. Therefore from the
definition of $L_f$ we have that

$$\tilde{L}_f(\overline\partial \phi_{\upharpoonright
W})=L_f(\overline\partial \phi_{\upharpoonright
W})=(-1)^{p+q+1}\,\int_{\text{Reg}\,\Omega} \phi\wedge f$$

\medskip
\noindent On the other hand from $(\ref{eq:Dfn})$ we have that

$$
\tilde{L}_f(\overline\partial \phi_{\upharpoonright
W})=\int_{\text{Reg}\,W} \overline\partial \phi \wedge
u=\int_{\text{Reg}\,\Omega} \overline\partial \phi \wedge u.
$$

\medskip
\noindent Putting the last two equalities together we obtain
$(\ref{eq:Lsolv})$.

\medskip
\noindent The last part of Theorem 5.3 follows directly from the
definition of $u$ and $(\ref{eq:L2normu'})$.

\medskip
\noindent

\medskip
\noindent{\bf{Remark:}} We can be more precise about the dependence
of $\text{supp}\,u$ and the constant $C$ (that appears in Theorem
5.3) on  $\text{supp}\,f$. Let $X,\Omega$ be as in theorem 5.3 and
let $N_0$ be a non-negative integer. There exists a positive integer
$N$ that depends on $N_0$ and $\Omega$ such that the following is
true: For every compact $K\subset \Omega$, there exists a compact
$K'\subset \Omega$ and a positive constant $C$ that depends on
$K,\,N,\,N_0,\,\Omega$ such that for every $(p,q)$ form $f$ with
$\text{supp}\,f\subset K$ and $\int_{K\setminus
A}|f|^2\,d_A^{N_0}\,dV<\infty$, $\overline\partial$-closed on
$\text{Reg}\,\Omega$, there exists a solution $u$ to
$\overline\partial u=f$ on $\text{Reg}\,\Omega$ with
$\text{supp}\,u\subset K'$ satisfying

$$
\int_{K'\setminus A} |u|^2\,d_A^{N}\,dV\le C\,\int_{K\setminus A}
|f|^2\,d_A^{N_0}\,dV.
$$

\medskip
\noindent The proof follows along the same lines as the proof of
Theorem 5.3, by taking as $p(v):=\left(\int_{K\setminus A}
|v|^2\,d_A^{-N_1}\,dV\right)^{\frac{1}{2}}$. Then there exist an
open, relatively compact subset $W$ of $\Omega$ (that depends on
$K$) a seminorm $q(w):=\int_{W\setminus A} |w|^2\,d_A^{-N}\,dV$ and
a positive constant $C'$ (that depends on $K$) such that the
equation $\overline\partial v=w$ has a solution $v$ satisfying
$p(v)\le C'\,q(w)$. The rest of the proof of Theorem 5.3 carries
over. The compact $K'$ is chosen to be $K':=\text{closure}(W)$ and
the constant $C=C'\,C_0$ where $C_0$ depends on $N,\,N_0,\,\Omega$.

\end{document}